\documentclass[11pt]{article}
   \usepackage{amsfonts}
\usepackage{latexsym,tikz}
\usepackage{amsmath}
 
\title{Multiple Dedekind Zeta Functions}
\author{Ivan Emilov Horozov}
\date{February 2, 2008}

\newcommand \nc {\newcommand}
\nc \proof {\noindent {\em{Proof.\/ }}} \nc \qed {$\Box$\hfill}
\newtheorem{theorem}{Theorem}[section]
\newtheorem{lemma}[theorem]{Lemma}
\newtheorem{proposition}[theorem]{Proposition}
\newtheorem{corollary}[theorem]{Corollary}
\newtheorem{definition}[theorem]{Definition}
\newtheorem{example}[theorem]{Example}
\newtheorem{remark}[theorem]{Remark}
\newtheorem{conjecture}[theorem]{Conjecture}
\newtheorem{question}[theorem]{Question}
\nc \bth[1] {\begin{theorem}\label{t#1} } \nc \ble[1]
{\begin{lemma}\label{l#1} } \nc \bpr[1]
{\begin{proposition}\label{p#1} } \nc \bco[1]
{\begin{corollary}\label{c#1} } \nc \bde[1]
{\begin{definition}\label{d#1}\rm } \nc \bex[1]
{\begin{example}\label{e#1}\rm } \nc \bre[1]
{\begin{remark}\label{r#1}\rm } \nc \bcon[1]
{\begin{conjecture}\label{con#1}\rm } \nc \bque[1]
{\begin{question}\label{que#1}\rm }
\nc {\eth} { \end{theorem} } \nc {\ele} { \end{lemma} } \nc
{\epr}{ \end{proposition} } \nc {\eco} { \end{corollary} } \nc
{\ede} {\end{definition} } \nc {\eex} { \end{example} } \nc {\ere}
{\end{remark} } \nc {\econ} { \end{conjecture} } \nc {\eque}
{\end{question} }

\def \A {{\mathbb A}}

\def \N {{\mathbb N}}
\def \Z {{\mathbb Z}}
\def \Q {{\mathbb Q}}
\def \R {{\mathbb R}}

\begin{document}

\title{{\LARGE\bf{Multiple Zeta Values and Ideles}}}

\author{
I. ~Horozov
\thanks{E-mail: ihorozov@math.wustl.edu}
\\ \hfill\\ \normalsize \textit{Department of Mathematics,}\\
\normalsize \textit{Washington University in St. Louis,}\\
}
\date{}
\maketitle
\begin{abstract}
In this paper we give two idelic representations of the multiple
zeta values - one using iterated integrals over the finite ideles and the other using iterated integrals over the idele class group. Each of the representations leads to a shuffle relation. Thus, we recover in a unified way the two types of shuffle relations of multiple zeta values via the iterated integrals over finite ideles  and via iterated integrals over the idele class group. 
\end{abstract}

\tableofcontents
\setcounter{section}{-1}
\section{Introduction}
Multiple zeta values are the values of a multiple zeta functions at the positive integers. 
Let us recall the Riemann zeta values: $$\zeta(k)=\sum_{n>0}\frac{1}{n^k}.$$ They were examined
first by Euler. 

Multiple zeta values of depth $d$ is
$$\zeta(k_1,\dots,k_d)=\sum_{0<n_1<\dots<n_d}
\frac{1}{n_1^{k_1}\dots n_d^{k_d}}.$$
These values were also examined by Euler.

Kontsevich expressed the
multiple zeta values as iterated integrals. 

\bth{0.1} (\cite{G})
Let
$k_1,\dots,k_d$ be $d$ positive integers with $k_d>1$. Then
$$\zeta(k_1,\dots,k_d)=
\int_0^1(\dots (\int_0^{x_3}(\int_0^{x_2}\frac{dx_1}{1-x_1})
\underbrace{\frac{dx_2}{x_2})\dots \frac{dx_{k_1}}{x_{k_1}}}_{k_1-1})
\frac{dx_{k_1+1}}{1-x_{k_1+1}})
\dots \frac{dx_{k_1+\dots+k_d}}{x_{k_1+\dots+k_d}}.$$
\eth
The iterated integral on the right hand side can be written as
$$\int\dots\int_{0<x_1<x_2<\dots<x_{k_1+\dots+k_d}<1}
\frac{dx_1}{1-x_1}\wedge
\underbrace{\frac{dx_2}{x_2}\wedge\dots\wedge \frac{dx_{k_1}}{x_{k_1}}}_{k_1-1}\wedge\frac{dx_{k_1+1}}{1-x_{k_1+1}}
\dots \wedge \frac{dx_{k_1+\dots+k_d}}{x_{k_1+\dots+k_d}}.$$

In this paper we give idelic interpretation of multiple zeta values (MZVs). 
The advantage of this approach is that the double shuffle relations for MZVs follow directly from the two idelic representations. One representation of a MZV is as an iterated integral over the finite ideles, which we define in this paper. The stuffle relations for MZV follows directly from that representation. The other representation of MZV is as an iterated integral over the idele class group. The shuffle relations for MZV follows directly from this representation. We prove the following:

\begin{theorem}
\label{thm main}
(Double shuffle relations via ideles) 

(a) $\zeta(k_1,\dots,k_d)=I(k_1,\dots,k_d)$,
where $I(k_1,\dots,k_d)$ is an iterated integral over the finite ideles $\A^\times_f$. Moreover, this this representation proves the stuffle relations among multiple zeta values.

(b) $\zeta(k_1,\dots,k_d)=J(k_1,\dots,k_d)$,
where $J(k_1,\dots,k_d)$ is an iterated integral over the ideles $\A^\times$. Moreover, this this representation proves the shuffle relations among multiple zeta values.

\end{theorem}

In Section 1, we recall particular Haar measures over local fields. We relate them to the Riemann zeta function following Tate \cite{T}. Then we define iterated integrals over the finite ideles and prove a stuffle relation formula for them. We also define iterated integrals over the idele class group and prove a shuffle relation for them.

In Sections 2, we express any MZV as an iterated integral over the finite ideles by applying the results from Section 1. Similarly, we express any MZV as an iterated integral over the idele class group. These two representations together with the stuffle and the shuffle relations from Section 1 allows us to prove the double shuffle relations for MZVs, Theorem \ref{thm main}, using ideles.

This approach has many useful features. First, it generalizes Tate's results to definition multiple zeta functions via idelic integration (Theorem 2.1). Second, we expect that similar results to hold for function fields over a finite field. That could lead to double shuffle relations for analogues of MZV over function fields over a finite field.
Finally, we expect similar idelic irepresentations hold for multiple Dedekind zeta values, (defined in \cite{MDZF}). They should lead to the two types of shuffle relations among multiple Dedkind zeta values \cite{shuffle}.

{\bf Acknowledgment:} I would
like to thank also to Professor Goncharov for the interest in
this work and for the encouragement he gave me.

This work was initiated at Max-Planck Institute f\"ur Mathematik.
I am very grateful for the stimulating atmosphere, created there.
Many thanks
are due to the University of Durham for the kind hospitality
during the academic year 2005-2006, when part of this work was
done, and to the Arithmetic Algebraic Geometry Marie Curie Network
for the financial support.

\section{Iterated integrals over the finite ideles or over the idele class group}
According to Tate's thesis \cite{T} the Riemann zeta function can be written as
a product of $p$-adic integrals. We want to expand this representation to
 multiple zeta functions as iterated adelic integrals. More precisely, iterated integrals over the finite ideles and over the idele class group.

Let $x_p$ be a $p$-adic number.
Let $|x_p|_p$ be the normalized $p$-adic norm so that $|p|_p=p^{-1}$.
Let $d_px_p$ be the additive $p$-adic Haar measure so that
$$\int_{\Z_p}d_px_p=1.$$
Let $d_p^\times x_p$ be the Haar measure on $\Q_p^\times$
normalized so that
$$\int_{\Z_p^{\times}}d_p^{\times} x_p = 1.$$
The relation between the two measures is the following
$$d_p^{\times} x_p = \frac{p-1}{p}\frac{dx_p}{|x_p|_p}.$$
Indeed, $\frac{dx_p}{|x_p|_p}$ is a multiplicative Haar measure.
Also, $$\Z_p-\{0\}=\bigcup_{k=0}^\infty p^k\Z_p^\times,$$ and
$$\int_{p^k\Z_p^\times}dx_p=p^{-k}\int_{\Z_p^\times}dx_p.$$ We have,
$$1=\int_{\Z_p}dx=\sum_{k=0}^\infty\int_{p^k\Z_p^\times}dx_p=
\sum_{k=0}^\infty p^{-k}\int_{\Z_p^\times}dx_p=\frac{p}{p-1}\int_{\Z_p^\times}dx_p.$$
Then
$$\int_{\Z_p^\times}dx_p=\frac{p-1}{p}.$$
Therefore,
$$\int_{\Z_p^\times}\frac{dx_p}{|x_p|_p}=\int_{\Z_p^\times}dx_p=\frac{p-1}{p}
=\frac{p-1}{p}\int_{\Z_p^{\times}}d_p^{\times} x_p.$$

Let $E_p(x_p)$ be a function defined on $\Q_p^{\times}$ by
$$E_p(x_p)= \left\{\begin{tabular}{ll}
$1$  & $x\in \Z_p-\{0\}$, \\
$0$     & otherwise \\
\end{tabular}\right.
    $$
The local factor of the Riemann zeta function is given by
$$\frac{1}{1-p^{-s}}=\int_{\Q_p^\times}E_p(x_p)|x_p|^s_p d_p^{\times}x_p .$$
Indeed, for fixed value of $k$ the integrant $E_p(x_p)|x_p|^s_p$ is constant
on the set $p^k\Z_p^\times$. For $k<0$ we have $E_p(x_p)=0$. For $k\geq 0$ we have
$$\int_{p^k\Z_p^\times}E_p(x_p)|x_p|^s_pd_p^\times x_p=\int_{p^k\Z_p^\times}|x_p|^s_pd_p^\times x_p=p^{-ks}.$$
Also, $\Z_p-\{0\}=\bigcup_{k=0}^\infty p^k \Z_p^\times$.
$$\int_{\Q_p^\times}E_p(x_p)|x_p|^s_p d_p^{\times}x_p=
\int_{\Z_p-\{0\}}|x_p|^s_p d_p^{\times}x_p=
\sum_{k=0}^\infty \int_{p^k\Z_p^\times}|x_p|^s_p d_p^{\times}x_p=
\sum_{k=0}^\infty p^{-ks}=\frac{1}{1-p^{-s}}$$

Denote by $x_\infty$ an element of $\R$. Let $|x_\infty|_\infty$ be the norm which
by definition is the absolute value of the real number. (We save the notation
$|x|$ for a norm of an idele.) Let $dx_\infty$ be the Haar measure on the additive group
of the real numbers.
Consider the multiplicative Haar measure on
$\R^\times$, namely,
$$\frac{dx_\infty}{|x_\infty|_\infty}$$
Let 
\begin{equation}
E_{\infty}(x_\infty)=\left\{
\begin{tabular}{ll}
$e^{-x_\infty}$ & for $x_\infty>0$,\\
$0$ & for $x_\infty<0$.
\end{tabular}\right.
\end{equation}
We are going to integrate
$E_\infty(x_\infty)$ with respect to the multiplicative measure.
The Mellin transform of
$E_\infty(x_\infty)$
gives the Gamma function
$$\int_{\R-\{0\}}E_{\infty}(x)|x|^s\frac{dx}{|x|}=\Gamma(s).$$


Let $x\in \hat{\Z}$. Denote by $|x|_f$ the product of all the $p$-adic norms.
Namely,  $$|x|_f=\prod_{p:finite}|x|_p.$$ Let also
$$E_f(x)=\prod_{p:finite}E_p(x_p).$$ Denote by
$$d_f^\times x_f=\prod_p d_p^\times x_p,$$
the multiplicative measure on the finite ideles given as product of all
local multiplicative measures over $\Q_p$ for all primes $p$.

\begin{definition}
\label{def it int over finite ideles}
Let $G_1,\dots,G_n$ be integrable functions on the finite ideles $\A_f^\times$ with support on $\hat{\Z}-\{0\}$. Assume that the function $G_i$ is constant on each set with a fixed norm, that is
$G_i$ is constant on the set $\{x_f\in \A_f^\times \,|\, |x_f|_f=q\}$ for a chosen $q\in \Q$.
We define an iterated integral over the finite ideles in the following way
\[\int_{\A_f^\times}G_1\circ\dots\circ G_n=\int_{|x_1|_f>\dots>|x_n|_f>0}G_1(x_1)\dots G_n(x_n)dx_1 dx_2\dots dx_n.\]
\end{definition}

\begin{definition}
\label{def stuffle}
We define the set $St(i,j)$ of stuffles $\sigma$ between an ordered  set with $i$ elements of $\Z$
\[0<k_1<\dots<k_i\]
and an ordered set with $j$ elements of $\Z$
\[0<k_{i+1}<\dots<k_{i+j}\]
to be all possible choices of ordered sets with $l$ elements of $\Z$
\[0<m_1<\dots<m_l\]
such that 
\[\{m_1,\dots,m_l\}=\{k_1,\dots,k_i\}\cup\{k_{i+1},\dots,k_{i+j}\}\]
where the two sets on the right hand side might not be disjoint.
If $\sigma$ is in $St(i,j)$ we denote by $\sigma$ both the map of inclusion  
\[\{k_1,\dots,k_i\}\rightarrow \{m_1,\dots,m_l\}\]
and the map of inclusion
\[\{k_{i+1},\dots,k_{i+j}\}\rightarrow \{m_1,\dots,m_l\}.\]
\end{definition}

\begin{theorem}
\label{thm stuffle}
(Stuffle relations for iterated integrals over the finite ideles)
Let $G_1,\dots,G_{i+j}$ be integrable functions on the finite ideles such that each of them is constant on each set of fixed norm, that is $G_i$ is constant on the set $\{x_f\in \A_f^\times \,|\, |x_f|_f=q\}$ for a chosen $q\in \Q$.
Then
\[\int_{\A_f^\times}G_1\circ\dots\circ G_i\times \int_{\A_f^\times}G_{i+1}\circ\dots\circ G_{i+j}
=\sum_{\sigma\in St(i,j)}G^\sigma_1\circ\dots\circ G^\sigma_l\]
where the sum is over all the stuffles $\sigma\in St(i,j)$ and 
\[G^\sigma_n=
\left\{
\begin{tabular}{llll}
$G_s$ & if $\sigma(k_s)=m_n$ for only one value of $s$;\\
$G_sG_t$ & if $\sigma(k_s)=m_n$ and $\sigma(k_t)=m_n$ for different values of $s$ and $t$.
\end{tabular}
\right.
\]
\end{theorem}
\proof Define $g_i(k)$ to be the integral of $G_i$ over the subset of the finite ideles of norm $1/k$, that is,  
\[g_i(k)=\int_{|x|_f=\frac{1}{k}} G_i(x)dx.\]
Since each of the functions $G_i$ is constant on subsets of the finite ideles with fixed norm, we obtain that
\[\int_{A^\times_f}G_1\circ\dots\circ G_i=\sum_{0<k_1<\dots<k_i}g_1(k_1)\dots g_i(k_i).\]
Then the stuffle relations among iterated integrals over the finite ideles are reduces to the stuffle relations of infinite series, namely,
\begin{align*}
\sum_{0<k_1<\dots<k_i}g_1(k_1)\dots g_i(k_i)
\sum_{0<k_{i+1}<\dots<k_{i+j}}g_{i+1}(k_{i+1})\dots g_i(k_{i+j})
=\\
=
\sum_{\sigma\in St(i,j)} 
\sum_{0<m_1<\dots <m_l}g^\sigma_1(m_1)\dots g^\sigma_l(m_l),
\end{align*}
where the first sum is over all the stuffles $\sigma\in St(i,j)$ and 
\[g^\sigma_n=
\left\{
\begin{tabular}{llll}
$g_s$ & if $\sigma(k_s)=m_n$ for only one value of $s$;\\
$g_sg_t$ & if $\sigma(k_s)=m_n$ and $\sigma(k_t)=m_n$ for different values of $s$ and $t$.
\end{tabular}
\right.
\]
\qed

\begin{definition}
\label{def it int over idele class group}
Let $G_1,\dots,G_n$ be integrable functions on the idele classes $\A^\times/\Q^\times$.
Assume that the function $G_i$ is constant on each set with a fixed norm, that is,
$G_i$ is constant on $\{x\in \A^\times \,|\, |x|=r\}$ for a chosen $r\in \R_{>0}$.
We define an iterated integral over the ideles in the following way
\[\int_{\A^\times/ \Q^{\times}}G_1\circ\dots\circ G_n=\int_{|x_1|>\dots>|x_n|>0}G_1(x_1)\dots G_n(x_n)dx_1 dx_2\dots dx_n.\]
\end{definition}

\begin{definition}
\label{def shuffle}
We define a shuffle of two ordered sets 
$\{1,\dots,i\}$ and $\{i+1,\dots,i+j\}$ as a permutation $\sigma$ of $i+j$ elements such that $\sigma$ satisfies the conditions
\[\sigma(1)<\dots<\sigma(i)\]
and
\[\sigma(i+1)<\dots<\sigma(i+j).\]
We denote by $Sh(i,j)$ the set of all shuffles between an ordered set of $i$ elements and an ordered set of $J$ elements.
\end{definition}
The key difference between a shuffle and a stuffle is that in the shuffle we have two disjoint sets while in the stuffle the two sets might have common elements. 

\begin{theorem}
\label{thm shuffle}
(Shuffle relations for iterated integrals over the idele class group)
Let $G_1,\dots,G_{i+j}$ be integrable functions on the finite ideles such that each of them is constant on each set of fixed norm, that is $G_i$ is constant on $\{x_f\in \A_f^\times \,|\, |x_f|_f=q\}$ for a chosen $q\in \Q$.
Then
\[\int_{\A_f^\times}G_1\circ\dots\circ G_i\times \int_{\A_f^\times}G_{i+1}\circ\dots\circ G_{i+j}
=\sum_{\sigma\in Sh(i,j)}G_{\sigma(1)}\circ\dots\circ G_{\sigma(i+j)}\]
where the sum is over all the shuffles $\sigma\in Sh(i,j)$.
\end{theorem}
\proof Let $\A^\times/\Q^\times\rightarrow \R_{>0}$ be the map under the norm of an idele. Since $G_l$ for $l=1,\dots,i+j$ are constant on subsets of the idele class group of fixed norm, we can define
\[g_i(r)=G_i(x)\] for $|x|=r$.
Since each of the functions $G_i$ are constant on subsets of the finite ideles with constant norm, we obtain that
\[\int_{A^\times_f}G_1\circ\dots\circ G_i=\sum_{0<r_1<\dots r_i}g_1(r_1)\frac{dr_1}{r_1}\dots g_i(r_i)\frac{dr_i}{r_i}.\]
Then the shuffle relation among iterated integrals over the idele class group reduces to the shuffle relation of iterated path integrals
\begin{align*}
&\int_{k_1>\dots>k_i>0}g_1(r_1)\dots g_i(r_i)\frac{dr_1}{r_1}\dots\frac{dr_i}{r_i}\times\\
&\times
\int_{r_{i+1}>\dots>r_{i+j}>0}g_{i+1}(r_{i+1})\dots g_i(r_{i+j})\frac{dr_{i+1}}{r_{i+1}}\dots\frac{dr_{i+j}}{r_{i+j}}
=\\
=&
\sum_{\sigma\in Sh(i,j)} 
\int_{r_1>\dots >r_{i+j}>0}g_{\sigma(1)}(r_1)\dots g_{\sigma(i+j)}(r_{i+j})
\frac{dr_1}{r_1}\dots\frac{dr_{i+j}}{r_{i+j}}.
\end{align*}
\qed

\section{Double shuffle relations for multiple zeta values via ideles}
\label{subsection MZV and ideles}

Let us iterate the function $E_f(x_f)|x_f|^s$ over the finite ideles.
We define 
\[I(s_1,s_2)=\int_{|x_1|_f>|x_2|_f}E_f(x_1)|x_1|_f^{s_1}d_f^\times x_1
E_f(x_2)|x_2|_f^{s_2}d_f^\times x_2.\]
If we iterate $d$ times, we obtain a multiple zeta function of depth $d$. Namely,
if we set
\begin{equation}
\label{eq finite adeles}
I(s_1,\dots ,s_d)=\int_{|x_1|_f>\dots>|x_d|_f}E_f(x_1)|x_1|_f^{s_1}d_f^\times x_1
\dots E_f(x_d)|x_d|_f^{s_d}d_f^\times x_d
\end{equation}

Recall the definition of multiple zeta function,
\[\zeta(s_1,\dots,s_d)=\sum_{0<n_1<\dots<n_d}\frac{1}{n_1^{s_1}\dots n_d^{s_d}}.\]
\begin{theorem}
\label{thm I}
Multiple zeta functions can be represented as iterated integrals over the finite ideles, namely,
\[\zeta(s_1,\dots,s_d)=I(s_1,\dots,s_d).\]
\end{theorem}
\proof
We are going to prove the theorem for $d=2$. For larger values of $d$ the proof is essentially the same.

Restrict the integral to the domain where the support of $f_f$ is not zero.
For any such idele we have $|x|_f=1/n$ for some $n\in \N$. Therefore,
$$\begin{tabular}{l}
$I(s_1,s_2)=$\\
\\
$=\int_{|x_1|_f>|x_2|_f}E_f(x_1)|x_1|_f^{s_1}d_f^\times x_1
E_f(x_2)|x_2|_f^{s_2}d_f^\times x_2=$\\
\\
$=\sum_{0<n_1<n_2}\int_{n_1\hat\Z^\times \times n_2\hat\Z^\times}
|x_1|_f^{s_1}d_f^\times x_1
|x_2|_f^{s_2}d_f^\times x_2=$\\
\\
$=\sum_{0<n_1<n_2}\frac{1}{n_1^{s_1}n_2^{s_2}}=$\\
\\
$=\zeta(s_1,s_2)$
\end{tabular}
$$
is a double zeta function. \qed


Let $x$ be an element of the adeles over $\Q$.
We are going to write $x_\infty$ for the infinite coordinate
of the adele $x$, and $x_p$ for the $p$-adic coordinate. Consider the function
$$E(x)=E_\infty(x_\infty)\prod_p E_p(x_p)$$
Let $d^{\times}x$ be a multiplicative
measure on the ideles given by the product
of local multiplicative measures considered above for all of the local fields.
Let $\A$ be the adels over the rational numbers $\Q$.
We are going to integrate over the idele class group $\A^\times/\Q^\times$.
For an idele $x \in \A^\times$ let
\[|x|=|x_\infty|_\infty\prod_p |x_p|_p\]
be the product of all the local valuations and let
\[e(x)=\sum_{q\in \Q^\times}E(qx).\] Define also
\[w_{1}(x)=e(x)d^\times x\]
and let
\[w_{0}(x)= d^\times x\]
be measures on $\A^\times$.

We define the following integral
\[J(m,n)=\int_{|x_1|>|x_2|>\dots|x_{m+n}|}
w_1(x_1)w_0(x_2)\dots w_0(x_m)w_{1}(x_{m+1})w_0(x_{m+2})\dots w_0(x_{m+n})
\]
More generally, we define
\begin{align}
\label{eq adeles}
&J(n_1,\dots , n_d)=\\
\nonumber
&=\int_{|x_1|>|x_2|>\dots>|x_{n_1+\dots+n_d}|}
\prod_{i=1}^d \left(w_1(x_{n_1+\dots+n_{i-1}+1})\left(\prod_{j=2}^{n_i}w_0(x_{n_1+\dots+n_{i-1}+j})\right)\right).
\end{align}

\begin{theorem}
\label{thm J} The multiple zeta values can be represented as iterated integrals over the idele class group, namely,
\[\zeta(n_1,\dots,n_d)=J(n_1,\dots,n_d).\]
\end{theorem}

\proof We are going to prove the Theorem for $d=2$. For larger $d$ the proof is essentially the same.

We want to modify the function $E(x)$ in the definition of
$J(m,n)$ so that the two new functions are defined over
$\A^\times/\Q^\times$ and over $\R^\times_{>0}$, respectively.

Denote by $\bar{x}$ the projection of an idele $x$ to an element of $\A^\times/\Q^\times$.
Let $$e(\bar{x})=\sum_{q\in \Q^\times}E(qx).$$ Note that $|qx|=|x|$ for $q\in \Q$ and
$x\in \A$. Denote by $|\bar{x}|:=|x|$ the norm of an element in $\A^\times/\Q^\times$.
 Recall that $\Q^\times$ is a discrete subgroup of $\A^\times$.
For that reason we can take the same measure $d^\times x$ on the set
$\A^\times/\Q^\times$. 

Now we define the corresponding function $c(t)$ for $t\in \R$ and $t>0$.
The function $e$ is constant on each set $(qt,q\hat{\Z}^\times)$,
where $t\in \R^\times_{>0}$ is fixed, $\hat{\Z}=\prod_p \Z_p$ and $q$ varies in $\Q^\times$. We can put $c(t)$ to be the value of
$e$ on any of the elements in the set $(qt,q\hat{\Z}^\times)$.
Let $x$ be an idele and $t=|x|$ be positive real number.
Let us examine more carefully the relation between $e(x)$ and
$c(t)$. For some $q\in \Q^\times$ and
$t\in \R^\times$ we have $x\in (qt,q\hat{\Z}^\times)$. If $q$ is
not an integer then $e(x)=0$. Also, $E_\infty(\mbox{negative }t)=0$. For these
reasons we can sum over all positive integers. Let $x_f$ be the
finite idele of $x$. That is, $x_f$ consists of
 all coordinates of $x$ except the coordinate corresponding to the infinite place. Let, also
$$E_f(x_f)=\prod_p E_p(x_p),$$ where the product is over all primes (finite places) $p$. Denote by $d^\times x_f$ the product
of multiplicative Haar measure of $\Q_p^\times$ over all primes $p$.
Denote, also, by $\hat{\Z}$ the product of the $p$-adic integers over all
primes $p$.
Then $$c(|x|)=\sum_{n\in \N}E_\infty(nr)\times
\int_{n\hat\Z^\times}E_f(x_f)d^\times x_f.$$
For $n\in \N$ the integral becomes
$$\int_{n\hat\Z^\times}E_f(x_f)d^\times x_f=1.$$
For the infinite place we have
$$E_\infty(nt)=e^{-nt}.$$
Therefore,
$$c(t)=\sum_{n=1}^\infty e^{-nt }.$$
Iteration of the measures $w_1$'s and $w_0$'s can be written as iteration of
$c(t)$:
\begin{align}
\nonumber
J(m,n)&
=\int_{|x_1|>|x_2|>\dots|x_{m+n}|}
w_1(x_1)w_0(x_2)\dots w_0(x_m)w_{1}(x_{m+1})\dots w_0(x_{m+n)}=\\
&=\int_{t_1>t_2>\dots>t_{m+n}>0}c(t_1)c(t_{m+1})dt_1\dots dt_{m+n}=
\zeta(m,n)
\end{align}
\qed

\proof (of the main Theorem \ref{thm main})
Part (a) follows from the explicit formulas for MZV in terms of the finite ideles from Theorem \ref{thm I} and from the corresponding stuffle relations from Theorem \ref{thm stuffle}. Part (b) follows from the explicit formulas for MZV in terms of idele class group from Theorem \ref{thm J} and from the corresponding shuffle relations from Theorem \ref{thm shuffle}.
\qed

\renewcommand{\em}{\textrm}
\begin{small}

\end{small}

\begin{thebibliography}{BHY1}

   \bibitem[L]{L}
Lang S.: {\em{Algebraic Number Theory,} 2nd ed. New York, Springer-Verlag, 1994.}
 
 \bibitem[G]{G}
Goncharov, A. B.:
{\em{Multiple polylogarithms and mixed Tate motives,}
math.AG/0103059, 82 pages.}


\bibitem[H1]{adeles}
Horozov, I. : {\em{Multiple Zeta Functions, Modular Forms and Adeles},
arXiv:math/0611849 [math.NT], 12 pages.}


\bibitem[H2]{MDZF}
Horozov, I.: {\em{Multiple Dedekind Zeta Functions},
arXiv:1101.1594 [math.NT], Journal f\"ur die reine und angewandte Mathematik (Crelle's Journal), to appear}.
	
\bibitem[H3]{shuffle}
Horozov, I: {\em{Double Shuffle Relations for Multiple Dedekind Zeta Values},
arXiv:1311.4019 [math.NT], 30 pages, submitted.}



    \bibitem[T]{T}
Tate, J.: {\em{Fourier Analysis in Number fields and Hecke's Zeta-Functions,}
Algebraic Number Theory, J.W.S. Cassels, A. Fr\"ohlich., Thompson Book Company
inc. Washington D.C., 1967, p. 305-347.}


    \end{thebibliography}
\end{document}